[1994/12/01]% LaTeX date must December 1994 or later
\documentclass[11pt,reqno]{article}

\title{Recurrence and Transience  for Branching Random Walks in an iid
Random Environment}
\author{Sebastian M\"uller}
\date{15.05.06}
\usepackage{amsmath,amsthm,amsfonts}
\usepackage{t1enc}
\usepackage[ps2pdf]{hyperref}
\chardef\bslash=`\\ % p. 424, TeXbook
\newtheorem{thm}{Theorem}[section]
\newtheorem{cor}[thm]{Corollary}
\newtheorem{lem}[thm]{Lemma}

\theoremstyle{definition}
\newtheorem{defn}{Definition}[section]

\theoremstyle{remark}
\newtheorem{rem}{Remark}[section]
\theoremstyle{example}
\newtheorem{ex}{Example}[section]

\newcommand{\eps}{\varepsilon}

\newcommand{\R}{\mathbb{R}}
\newcommand{\Z}{\mathbb{Z}}
\newcommand{\N}{\mathbb{N}}
\renewcommand{\P}{\mathbb{P}}
\newcommand{\E}{\mathbb{E}}

\newcommand{\eval}[2][\right]{\relax
  \ifx#1\right\relax \left.\fi#2#1\rvert}

\begin{document}
\maketitle {\abstract We give three different criteria for
transience of a Branching Markov Chain. These conditions enable us
to give a classification of Branching Random Walks in Random
Environment (BRWRE) on Cayley Graphs in recurrence and transience.
This classification is stated explicitly for BRWRE on $\Z^d.$
Furthermore, we emphasize  the interplay between Branching Markov
Chains and the spectral radius. We prove properties of the
spectral radius of the Random Walk in Random Environment with the
help of appropriate Branching Markov Chains.}
\renewcommand{\sectionmark}[1]{}
\section{Introduction}

A Branching Markov Chain (BMC)  is  a system of particles in
discrete time. The BMC starts with one particle in an arbitrary
starting position $x_s.$ At each time  particles split up in
offspring particles independently according to some probability
distributions $\mu$, that may depend on the locations of the
particles. The new particles then move independently according to
a Markov Chain  (MC).

An irreducible MC is either recurrent or transient:  either all or
none states are visited infinitely often.  It turns out that this
dichotomy breaks down for BMC.  Let $\alpha(x)$ be the probability
that, starting the BMC in $x_s=x$, the state $x$ is hit infinitely
often by some particles. There are three possible regimes:
\emph{transient} $(\alpha(x)=0~ \forall x)$, \emph{weakly
recurrent} $(0<\alpha(x)<1~ \forall x)$ and \emph{strongly
recurrent} $(\alpha(x)=1~ \forall x).$ For a discussion of these
results, a more detailed introduction and references we refer to
\cite{gantert2004}. In Theorem \ref{thm:1} we give equivalent
criteria for the transience of BMC. The interplay of these
different criteria is central in our development.

If the underlying MC is a Random Walk (RW) we speak of a Branching
Random Walk (BRW). A BRW on a Cayley Graph is either transient or
strongly recurrent if the offspring distribution is constant, i.e.
$\mu(x)=\mu$  for all vertices $x,$ with mean offspring
$m=\sum_{k\geq 1} k \mu_k,$ see \cite{gantert2004}. In particular,
we have that $m\leq 1/\rho$ implies transience and $m>1/\rho$
implies strong recurrence, where $\rho$ is the spectral radius of
the RW. We show that these results can be generalized to Branching
Random Walk in Random Environment (BRWRE). We consider iid
environments and assume the branching and the transition mechanism
to be independent, too. We obtain a necessary and sufficient
condition for  transience, see Theorem \ref{thm:5}. This condition
depends only on some extremal points of the support of the
environment. In particular, we give an explicit criterion for
transience and strong recurrence for BRWRE on $\Z^d$, see
Corollary \ref{cor:6}. This answers a question asked in
\cite{comets05}. We refer to \cite{comets05} for an investigation
of a more general model of BRWRE where the branching and movement
can be dependent. An additional purpose of this paper is to
emphasize the interplay between the behavior of the BRW and the
spectral radius of the underlying RW. On  one hand, the critical
mean offspring equals the inverse spectral radius of the RW. On
the other hand, we can use BRW to derive properties of the
spectral radius, see the proof of Lemma \ref{lem:5}.

\section{Preliminaries}

Let $G$ be  a finitely generated group. Unless $G$ is abelian, we
write the group operation multiplicatively. Let $S$ be a finite
symmetric set of generators of $G$ and $q$ a probability measure
on $S.$ The Cayley Graph $X(G,S)$ with respect to $S$ has vertex
set $G$, and two vertices $x,y\in G$ are connected if and only if
$x^{-1}y\in S.$ The Random Walk (RW) on $X(G,S)$ with transition
probabilities $q$ is the Markov Chain with state space $X=G$ and
transition probabilities
$$p(x,y)=q(x^{-1}y)\quad \mbox{for }x^{-1}y\in S$$ and $0$ otherwise.
 The $n-$step transition
probabilities are
$$p^{(n)}(x,y)=q^n(x^{-1}y),$$
where $q^n$ is the $n-$fold convolution of $q$ with itself. We
start the RW in a starting position $x_s.$

We introduce the Random Environment. Let $\mathcal{M}$ be the
collection of all probability measures on $S.$ Let
$(\omega_x)_{x\in X}$ be a collection of iid random variables with
values in $\mathcal{M}$ which serve as an environment. For
 each realization $\omega:=(\omega_x)_{x\in X}$ of this environment, we define a
Markov Chain $(X_n)_{n\in \N}$ on $X=G$ with $X_0=x_s$ and
$$\P_\omega(X_{n+1}=y| X_n=x):=p_\omega(x,y):=\omega_x(x^{-1} y)\quad
\forall n\geq 1.$$ We denote by $P_\omega$  the transition kernel
of the Markov chain on the state space $X.$

Let $\eta$ be the distribution of this environment. We assume that
$\eta$ is a product measure with  one-dimensional marginal $Q.$
The support of $Q$ is denoted by $\mathcal{K}$ and its convex hull
by $\hat{\mathcal{K}}.$ Throughout this note we assume the
following condition on $Q:$
\begin{eqnarray}\label{eq:cond}
      Q\{\omega: \omega(s)>\gamma~\forall
    s\in S'\}=1 \mbox{ for some } \gamma>0,
\end{eqnarray} where $S'\subseteq S$ is a minimal set of
generators. We assume it to ensure the irreducibility of a RW with
transition probabilities $q\in {\mathcal{K}}.$

We recall the definition of the  spectral  radius of an
irreducible Markov Chain $(X,P):$
\begin{equation}
\rho(P):=\limsup_{n\rightarrow\infty}
\left(p^{(n)}(x,x)\right)^{1/n} \in (0,1],
\end{equation}
where $p^{(n)}(x,x)$ is the probability to get from $x$ to $x$ in
$n$ steps. The following characterization of the spectral radius
in terms of $t-$superharmonic functions is crucial for our
classification:
\begin{lem}\label{lem:1}
$$\rho(P)=\min\{t>0:~\exists\, f(\cdot)>0\mbox{ such that }Pf\leq t f\}$$
\end{lem}
For the proof and more  information on RW on groups we refer to
\cite{woess}.

Due to the symmetry of the Cayley Graph and the independence of
the environment we find:

\begin{lem}\label{lem:2}We have $$
\rho:=\rho(P_\omega)=\sup_\sigma \rho(P_\sigma) ~for~
\eta\mbox{-a.a.}~\omega,$$ where the $\sup$ is over all possible
realizations $\sigma=(\sigma_x)_{x\in X}$ with
$\sigma_x\in\mathcal{K}.$
\end{lem}

\begin{proof}
It is obvious that $\rho(P_\omega)\leq\sup_\sigma\rho(P_\sigma).$
For the opposite inequality  let $\eps>0$ and $\tau$ a possible
realization such that $\rho(P_\tau) \geq \sup_\sigma
\rho(P_\sigma)-2\eps.$ Hence, for any $x\in X$ there exists
$n\in\N$ such that:
\begin{equation}\label{eq:lem:2}
\left(p_\tau^{(n)}(x,x)\right)^{1/n}\geq \sup_\sigma
\rho(P_\sigma)-\eps.
\end{equation}
We have $$ p_\tau^{(n)}(x,x)=\sum_{x=x_0,\ldots,x_n=x}
\prod_{i=0}^{n-1} p_\tau(x_i,x_{i+1}),$$ where the sum is over all
possible paths of length $n$ from $x$ to $x.$ Assume the
distribution $Q$ to be discrete, then we find, due to the
transitivity of the Cayley Graph and the independence of the
environment, for $\eta$-a.a. environments $\omega$ some other
vertex $y\in X$ such that
$$ p_\omega(y x_i ,y x_j )= p_\tau(x_i, x_j)\quad
\forall x_i,x_j \in \{z: d(x,z)\leq n/2\},$$ where
$d(\cdot,\cdot)$ is the usual graph distance. Hence,
$p_\omega^{(n)}(yx,yx)=p_\tau^{(n)}(x,x).$

Using the fact that $\rho(P_\omega)^n\geq p_\omega^{(n)}(z,z)$ for
all $z\in X$   we obtain with inequality (\ref{eq:lem:2}) that
$$\rho(P_\omega)\geq \left(p_\omega^{(n)}(yx,yx)\right)^{1/n} =
\left(p_\tau^{(n)}(x,x)\right)^{1/n}\geq \sup_\sigma
\rho(P_\sigma)-\eps \quad \forall \eps>0.$$

For the general case let $\delta>0,$ for $\eta$-a.a. environments
$\omega$ we find some $y\in X$ such that
$$ p_\omega(y x_i ,y x_j )\geq \frac1{1+\delta}~
p_\tau(x_i, x_j)\quad \forall x_i,x_j \in \{z: d(x,z)\leq n/2\}.$$
We have
\begin{eqnarray*}
  p_\omega^{(n)}( yx, yx)  &=& \sum_{ x=x_0,\ldots,x_n= x} \prod_{i=0}^{n-1}
p_\omega(y x_i ,y x_{i+1})\\
   &\geq&  \sum_{x=x_0,\ldots,x_n=x}
   \prod_{i=0}^{n-1}\frac1{1+\delta}~
p_\tau(x_i,x_{i+1}) \\
   &=& \left(\frac1{1+\delta}\right)^n p_\tau^{(n)}(x,x) \quad \forall
   \delta>0.
\end{eqnarray*}
Letting $\delta\to 0,$  this yields
$$\rho(P_\omega)\geq \sup_\sigma
\rho(P_\sigma)-\eps \quad \forall \eps>0.$$
\end{proof}
\begin{rem}
The fact that the spectral radius is constant for $\eta$-a.a.
realizations $\omega$ of the environment follows  directly from
the observation that
$$\rho(P_\omega)=\limsup_{n\rightarrow\infty}\left( p_\omega^{(n)}(x,x)\right)^{1/n}$$
does not depend on $x$ and hence, by ergodicity, is constant a.s..
\end{rem}

\subsection{Branching Markov Chains}\label{sect:BMC}
We introduce the model of Branching Markov Chain (BMC). Let
$(X,P)$ be an irreducible and infinite Markov Chain in discrete
time. For all $x\in X$ let
$$\mu(x)=\left(\mu_k(x)\right)_{k\geq 1}$$ be a sequence of non-negative numbers
satisfying
$$\sum_{k=1}^\infty \mu_k(x)=1 \mbox{ and } m(x):=\sum_{k=1}^\infty k \mu_k(x)<\infty.$$
We define the BMC $(X,P,\mu)$ with underlying Markov Chain $(X,P)$
and branching distribution $\mu=(\mu(x))_{x\in X}$ following
\cite{menshikov97}. At time $0$ we start with one particle in an
arbitrary starting position $x_s\in X.$  At time $1$ this particle
splits up in $k$ offspring particles  with probability
$\mu_k(x_s).$ Still at time $n=1,$ these $k$ offspring particles
then  move independently according to the Markov Chain $(X,P).$
The process is defined inductively. At each time each particle in
 position  $x$ splits up according to $\mu(x)$ and the offspring
particles move according to $(X,P).$ At any time, all particles
move and branch independently of the other particles and the
previous history of the process. Let $\eta(n)$ be the total number
of particles at time $n$ and let $x_i(n)$ denote the position of
the $i$th particle at time $n.$ Denote
$\P_x(\cdot)=\P(\cdot|x_s=x).$ We define  recurrence and
transience for  BMC as in \cite{gantert2004}:
\begin{defn} Let
\begin{equation}\label{eq:1} \alpha(x):=\P_x\left(
\sum_{n=1}^\infty
\sum_{i=1}^{\eta(n)}\mathbf{1}\{x_i(n)=x\}=\infty\right).
\end{equation}
 A BMC is {\it recurrent,} if $\alpha(x)>0$ for some  $x\in X,$
 {\it strongly recurrent,} if $\alpha(x)=1$ for some $x\in X$  and
{\it transient} otherwise.
\end{defn}
The definition does not depend on the starting position $x_s=x.$
In fact,
 $\alpha(x)>0,~ \alpha(x)=1$ and $\alpha(x)=0$ hold
either for all or none $x\in X,$ see  \cite{gantert2004}. We write
$\alpha>0$ if $\alpha(x)>0$ for all $x\in X$ and $\alpha\equiv 1$
and $\alpha\equiv 0$ respectively. In analogy to
\cite{menshikov97}, we introduce the following modified version of
the BMC. We fix some position $x_0\in X,$ which we denote the
origin of $X.$  The new process is like the original BMC at time
$n=1$ but is different for $n>1.$ After the first time step we
conceive the origin as {\it freezing}: if a particle reaches the
origin it stays there forever and stops splitting up. We denote
this new process with BMC*. The process BMC* is analogous to the
original process BMC except that $p(x_0,x_0)=1,~
p(x_0,x)=0~\forall x\neq x_0$ and $\mu_1(x_0)=1$ from the second
time step on. Let $\eta(n,x_0)$ be the number of particles at
position $x_0$ at time $n.$ We define the random variable
$\nu(x_0)$ as
$$\nu(x_0):=\lim_{n\rightarrow\infty} \eta(n,x_0).$$
The random variable $\nu(x_0)$ takes values in
$\{0,1,\ldots\}\cup\{\infty\}.$ We write $\E_x \nu(x_0)$ for the
expectation of $\nu(x_0)$ given that $x_s=x.$  Note that our
notation of $\eta(n,x_0)$ and $\nu(x_0)$ is different from the one
in \cite{gantert2004}. Since the choice of the origin may affect
the behavior of the BMC we keep track of the dependence of the
variables $\eta$ and $\nu$ on the choice of the origin and write
$\eta(n,x_0)$ and $\nu(x_0).$ Furthermore, our definition of the
process BMC* differs from the one given in \cite{menshikov97}. In
our definition the origin is not absorbing at time $n=1.$ These
modifications enable us to give the following three different
criteria for transience of BMC that hold for all irreducible and
infinite Markov Chains $(X,P).$

\begin{thm}\label{thm:1}
A BMC $(X,P,\mu)$ with $m(y)>1$ for some $y$ is transient if and
only if the three equivalent conditions hold:
\begin{enumerate}
     \item[(i)] $\E_x \nu(x)\leq 1$ for some/all $x\in X.$
     \item[(ii)] $\E_x \nu(x_0)<\infty$ for all $x,x_0\in X.$
     \item[(iii)] There
exists a strictly positive function $f(\cdot)$ such that
\begin{equation}\label{eq:1.1}
P f(x)\leq \frac{f(x)}{m(x)}\quad \forall x\in X.
\end{equation}

\end{enumerate}
\end{thm}
\begin{proof}
\underline{$(i)\Leftrightarrow ~\alpha\equiv 0:$} We start the BMC
in some $x_0\in X.$ The key idea of the proof is to observe that
the total number of particles ever returning to  $x_0$ can be
interpreted as the total number of progeny in a branching process
$(Z_n)_{n\geq 0}.$  Note that each particle has a unique ancestry
line which leads back to the starting particle at time $0$ at
$x_0.$ Let $Z_0:=1$ and $Z_1$ be the number of particles being the
first particle in their ancestry line to return to $x_0.$
Inductively we define $Z_n$ as the number of particles being the
$n$th particle in their ancestry line to return to $x_0.$ This
defines a Galton Watson process $(Z_n)$ with offspring
distribution $Z{\buildrel d \over =}Z_1.$ Observe that the
particles being the first in their ancestry line to visit $x_0$
are those that are frozen in $x_0$ in BMC* with origin $x_0.$
Hence, $Z_1{\buildrel d \over =}\nu(x_0)$ given that the process
starts in $x_0.$ Notice also that
$$ \alpha(x_0)=\P_{x_0}\left(
\sum_{n=1}^\infty
\sum_{i=1}^{\eta(n)}\mathbf{1}\{x_i(n)=x_0\}=\infty\right),$$ and
$$\sum_{n=1}^\infty Z_n =\sum_{n=1}^\infty \sum_{i=1}^{\eta(n,x_0)} \mathbf{1}\{x_i(n)=x_0\}.$$
If $\alpha(x_0)=0,$ then $\sum_{n=1}^\infty Z_n<\infty$ a.s.,
hence $(Z_n)$ is (sub-)critical and  $\E_{x_0} \nu(x_0)=E[Z]\leq
1.$

Now let $\E_{x_0}\nu(x_0)\leq 1.$ Since $m(y)>1$ for some $y$ and
hence $\P_{x_0} (\nu(x_0)>1)>0$ we have that $\E_{x_0}\nu(x_0)\leq
1$ implies that the process $(Z_n)$ dies out a.s.. Therefore,
$\alpha(x_0)=0.$ The claim follows since $\alpha(x)=0$ either
holds for all or none $x\in X.$

\underline{$(ii)\Leftrightarrow ~\alpha\equiv 0:$} Let first be
$\alpha\equiv 0$ and assume  that there exists $x_0$ and $x$ such
that $\E_x\nu(x_0)=\infty.$ Hence, $E_{x_0}\nu(x_0)=\infty,$ since
$(X,P)$ is irreducible. This contradicts $\alpha\equiv 0,$ since
$(i)$ is equivalent to $\alpha\equiv0.$

 In order to show the converse we use again a proof
by contradiction. We assume that $\alpha>0$ and show that
$\E_{x_0}\nu(x_0)=\infty$ for some $x_0.$  Due to (i) we have
$\E_x\nu(x)>1$ for some $x.$ Let $k$ be such that $\E_x
\eta(k,x)>1.$ Let us first assume that the Markov Chain $(X,P)$
has finite range, i.e. $|\{y:~ p(x,y)>0\}|<\infty$ for all $x\in
X.$   Since $(X,P)$ is infinite and has finite range we find some
$x_0$ such that $d(x,x_0)>k,$ where $d(x,x_0):=\inf\{n\in\N:
p^{(n)}(x,x_0)>0\}.$ We proceed as follows: we start a
 BMC* with origin $x_0$ in $x_0,$ with positive probability one
particle reaches $x,$ this particle initiates a supercritical
Galton-Watson process, $(\zeta_i)_{i\geq 0},$ of particles
visiting $x.$ Therefore, $x$ is visited infinitely often with
positive probability. The final step is then to show that this
implies that infinitely many particles are frozen in the origin
$x_0.$ This is clearly enough since $\P_{x_0}(
\eta(x_0)=\infty)>0$ implies that $\E_{x_0} \nu(x_0) =\infty.$

We define the Galton-Watson process ($\zeta_i)_{i\geq 0}.$ We
start a BMC* with origin $x_0$ with one particle in $x.$ Let
$\Psi_1$ be the particles being the first particles in their
ancestry line to return to $x$ before time $k.$ We define $\Psi_i$
inductively as the number of particles having an ancestor in
$\Psi_{i-1}$ and being  the first in the ancestry line of this
ancestor to return to $x$ in at most $k$ time steps. Clearly
$\zeta_0:=1$ and $\zeta_i:=|\Psi_i|,~i\geq 1,$ defines a
Galton-Watson process. Since $d(x,x_0)>k$ we have that
$$ E[\zeta_1]=\E_x\eta(k,x)>1.$$ Therefore, the process $(\zeta_i)$
is supercritical and survives with positive probability. This
implies that

$$\P_{x_0}\left( \sum_{n=1}^\infty
\sum_{i=1}^{\eta(n)}\mathbf{1}\{x_i(n)=x\}=\infty\right)>0.$$ We
shall prove $\P_{x_0}( \nu(x_0)=\infty)>0$ by showing that

\begin{equation}\label{eq:thm:1}\P_{x_0}\left( \sum_{n=1}^\infty
\sum_{i=1}^{\eta(n)}\mathbf{1}\{x_i(n)=x\}=\infty\mbox{ and }
\nu(x_0)<\infty\right)=0.\end{equation} We follow the line of the
proof of Lemma 3.3 in \cite{benjamini94}. Since $X$ is irreducible
we have $p^{(l)}(x,x_0)=\delta>0$ for some $l\in \N.$ Let $N,
M\in\N.$ The probability that there are times $M<n_1,\ldots,n_N$
with $n_{j-1}+l<n_j$ for all $1\leq j\leq N$ such that
$x_i(n_j)=x$ for some $1\leq i\leq \eta(n_j)$ for all $j$ but
$x_i(n)\neq x_0$ for all $n>M$ and all $1\leq i \leq \eta(n)$ is
at most $(1-\delta)^N.$ For this it is crucial that $x_0$ is the
only absorbing position in BMC*. Letting $N\rightarrow\infty,$
this yields
$$\P_{x_0}\left( \sum_{n=1}^\infty
\sum_{i=1}^{\eta(n)}\mathbf{1}\{x_i(n)=x\}=\infty\mbox{ and }
\nu(x_0)-\eta(M,x_0)=0\right)=0.$$ Let $A_M$ be the event in the
last formula. Notice that,
$$\bigcup_{M\geq 1} A_M =\left\{ \sum_{n=1}^\infty
\sum_{i=1}^{\eta(n)}\mathbf{1}\{x_i(n)=x\}=\infty\mbox{ and }
\nu(x_0)<\infty\right\}.$$ This gives equation (\ref{eq:thm:1}).

We now turn to the general case where $(X,P)$ is any irreducible
and infinite Markov Chain. We have used the finiteness of the
range of $(X,P)$ to ensure the existence of some $x_0$ such that
$d(x,x_0)>k.$ This was useful to easily bound  the mean offspring
of the process $(\zeta_i).$ For the general case we use
approximation arguments to show that $E[\zeta_1]>1$ for some
$x_0\in X.$ The  remaining part of the proof then follows the
arguments of the special case. Observe that a BMC and the
corresponding BMC* can also be considered as a Markov Chain on the
state space $X',$ namely the state space of all particle
configurations
$$c(n)=\{x_1(n), x_2(n),\ldots,x_{\eta(n)}(n)\}$$ where $x_i(n)$ is
the position of the the $i$th particle at time $n$ and $\eta(n)$
is the total number of particles at time $n.$ Let $C_k$ be the set
of all possible particle configurations of a BMC*, started in $x$
with origin $x,$ up to time $k.$ Let $c=(c(1),\ldots,c(k))\in C_k$
be a possible realization and let  $\eta_c$ denote the number of
particles frozen in $x_0$ at time $k$ for the realization $c.$
Observe now that
$$\E_x \eta(k,x)=\sum_{c\in C_k} \eta_c P'(c),$$
where $P'(c)$ is the probability that we see $c$ as realization of
the BMC* up to time $k.$ There exists a finite subset $C\subseteq
C_k$ s.t.
$$\sum_{c\in C} \eta_c P'(C)>1.$$
Since the number of different positions visited by some
realization $c\in C$ is finite we find some $x_0$ that is not
visited for all $c\in C.$ Hence, for this $x_0$ we have that
$E[\zeta_1]>1.$

\underline{$(iii)\Leftrightarrow ~\alpha\equiv 0:$} We first show
that $\alpha\equiv 0 $ implies $(iii).$  Let $x_0\in X.$ Due to
(i) and (ii) we have that  $\E_{x_0}\nu(x_0)\leq 1$  and $\E_{x}
\nu(x_0)<\infty$ for all $x.$ We show that $f(x):=\E_x\nu(x_0)>0$
satisfies inequality (\ref{eq:1.1}). For $x$ such that
$p(x,x_0)=0$ it is straightforward to show that even equality
holds in (\ref{eq:1.1}). Let $x$ such that $p(x,x_0)>0$ we have
\begin{eqnarray*}
 f(x)=\E_x\nu(x_0) &=& m(x) \left( \sum_{y\neq x_0}p(x,y) \E_y \nu(x_0)
+ p(x,x_0) \cdot 1 \right) \\
   &\geq& m(x)  \sum_y p(x,y) \E_y \nu(x_0) \\
   &=& m(x) Pf(x),
\end{eqnarray*}
 since $\E_{x_0} \nu(x_0)\leq 1.$

The proof that the existence of a function satisfying
(\ref{eq:1.1}) implies transience  is due to \cite{gantert2004}.
We give  a short
 sketch. Consider the  BMC* with origin $x_0$ and define
$$Q(n):=\sum_{i=1}^{\eta(n)} f(x_i(n)),$$
where $x_i(n)$ is the position of the $i$th particle at time $n.$
Observing  that $Q(n)$ is a positive  supermartingale that
converges a.s. to a random variable $Q_\infty$ and that
\begin{equation}\label{eq:q}
\nu(x_s)\leq \frac{Q_\infty}{f(x_s)}
\end{equation}
 we obtain
\begin{equation}\label{eq:1.11}
\E_{x_0}\nu(x_0)\leq \frac{\E_{x_0}Q_\infty}{f(x_0)}\leq
\frac{\E_{x_0}Q(0)}{f(x_0)}=\frac{f(x_0)}{f(x_0)}=
1.\end{equation}
\end{proof}

In particular if the mean offspring is constant, i.e. $m(x)=m
~\forall x\in X,$ we have, due to  Lemma \ref{lem:1}, the
following  result of \cite{gantert2004}:

\begin{thm}\label{thm:1a}
For a BMC with  underlying Markov chain $(X,P)$ and constant mean
offspring $m>1$, it holds that the BMC is transient if $m\leq
1/\rho(P)$ and recurrent if $m>1/\rho(P).$
\end{thm}

The next Theorem  follows from the argumentation of the proof of
  Theorem \ref{thm:1}, part (iii), and is due to
\cite{menshikov97}.
\begin{thm}\label{thm:4}
Let $x_0\in X.$ There exists a function $f>0$ satisfying
\begin{equation}\label{eq:thm4}
 P f(x)= \frac{f(x)}{m(x)}\quad \forall x\neq x_0
\end{equation}
if and only if $\E_{x} \nu(x_0)<\infty $ for all $x\neq x_0.$ In
this case a solution of (\ref{eq:thm4}) is given as $f(x):=\E_x
\nu(x_0)$ for $x\neq x_0$ and $f(x_0):=1.$
\end{thm}

 In order to transfer Theorem \ref{thm:1a} to  BMC
 with non-constant mean offspring we use coupling arguments.
 We couple a BMC $(X,P,\mu)$ with a suitable BMC $(X,P,\tilde\mu)$ with a
 given constant mean offspring $\widetilde m,$ $\widetilde m\geq m(x)~\forall x,$
 such that there are always everywhere more particles  in
$(X,P,\tilde\mu)$ than in $(X,P,\mu).$ We obtain that $(X,P,\mu)$
is transient if the coupled process $(X,P,\tilde\mu)$ is
transient. In order to describe the coupling we say that
$(X,P,\mu)$ consists of blue particles and $(X,P,\tilde\mu)$ of
blue and red particles. The coupling is defined such that the blue
particles in $(X,P,\tilde\mu)$ are a copy of the whole blue
process $(X,P,\mu).$ The red particles are considered as
supplementary particles.

 We choose the
distributions $(\tilde\mu(x))_{x\in X}$ as follows. For each $x\in
X$ let $l$ be the smallest integer such that $\mu_l(x)>0.$ Let
$\delta:=\widetilde{m}-m(x)$ and $n\in \N$ such that
$\mu_l(x)>\frac\delta{n}.$ Let
\begin{eqnarray*}
  \tilde\mu_l(x) &:=& \mu_l(x)-\frac\delta{n} \\
  \tilde\mu_{n+l}(x) &:=& \mu_{n+l}(x)+ \frac\delta{n}  \\
  \tilde\mu_i(x) &:=& \mu_i(x)\quad \forall i\notin \{l,n+l\}.
\end{eqnarray*} This defines a BMC $(X,P,\tilde\mu)$ with desired
mean offspring $\widetilde{m}(x)=\widetilde{m}~\forall x.$ We
couple the two processes inductively. Starting with one blue
particle in $x_s$ we produce $k$  blue offspring in $(X,P,\mu)$ if
$$ \sum_{i=1}^{k-1}\mu_i(x_s) \leq U < \sum_{i=1}^{k} \mu_i(x_s),$$
where $U$ is uniformly distributed on $[0,1].$ In
$(X,P,\tilde\mu)$ we produce $k,$ $k\neq l,$ blue  offspring  if $
\sum_{i=1}^{k-1}\mu_i(x_s) \leq U < \sum_{i=1}^{k} \mu_i(x_s),$
$l$ blue offspring if $U<\mu_l(x_s)-\frac\delta{n}$ and $l$ blue
and $n$ red offspring if $\mu_l(x_s)-\frac\delta{n}\leq U <
\mu_{l}(x_s).$ Note that in both processes we have the same number
of blue particles. The blue particles are coupled such that they
move to the same positions. For each of these blue particles we
start a new coupled branching mechanism as defined above
independent of each other and the previous history. The red
particles are not coupled to any particles in $(X,P,\mu)$ and
perform as a usual but red-colored BMC.

The same procedure can be applied to couple a BMC $(X,P,\mu)$ to a
suitable BMC $(X,P,\tilde\mu)$ with a given constant mean
offspring $\widetilde m,$ $\widetilde m\leq m(x)~\forall x,$ and
fewer particles than the original process. We  choose
$\tilde\mu(x)$ in the following way. For each $x$ let
$\delta:=m(x)-\widetilde m$ and $l$ be such that
$$\sum_{k=1}^l (k-1) \mu_k(x) \leq \delta < \sum_{k=1}^{l+1} (k-1)
\mu_k(x).$$ With $\gamma:=\delta - \sum_{k=1}^l (k-1) \mu_k(x)$ we
define the offspring distribution
\begin{eqnarray*}
  \tilde\mu_1(x) &:=& \mu_1(x)+\mu_2(x)+\cdots+\mu_l(x)+ \frac\gamma{l}\\
  \tilde\mu_i(x) &:=& 0\quad 2\leq i\leq l\\
  \tilde\mu_{l+1}(x) &:=& \mu_{l+1}(x) -\frac\gamma{l}\\
   \tilde\mu_{i}(x)&:=& \mu_i(x) \quad i>l+1
\end{eqnarray*} with desired mean offspring $\widetilde m(x)=\widetilde m.$

\begin{rem}\label{rem:couple}
A BMC  is transient if the coupled BMC with more particles is
transient. It is  recurrent if the coupled BMC with less particles
is recurrent.
\end{rem}

\section{BRWRE on Cayley Graphs}
In addition to the  environment that determines the random walk we
introduce a random environment determining the branching
mechanism. Let $\mathcal{B}$ be the set of all infinite positive
sequences $\mu=\left(\mu_k\right)_{k\geq 1}$ satisfying
$\sum_{k=1}^\infty \mu_k=1$ and $m(\mu):=\sum_{k=1}^\infty
k\mu_k<\infty.$ Let $\widetilde Q$ be a probability distribution
on $\mathcal{B}$ and set
\begin{equation}
\label{mstar} m^*:=\sup\{m(\mu):\mu\in {\rm supp}(\widetilde Q)\}.
 \end{equation}
Let $(\omega_x)_{x\in X}$ be a collection of iid random variables
with values in $\mathcal{M}$ and $(\mu_x)_{x\in X}$ be a
collection of iid random variables with values in $\mathcal{B}$
such that $(\omega_x)_{x\in X}$ and $(\mu_x)_{x\in X}$ are
independent, too. Let $\Theta$ be the corresponding product
measure with one-dimensional marginal $Q\times \widetilde Q.$ For
each realization $(\omega,\mu):=(\omega_x,\mu_x)_{x\in X}$ let
$P_\omega$ be the transition kernel of the underlying Markov Chain
and $\mu(x)=\mu_x.$  Thus, each realization $(\omega,\mu)$ defines
a BMC $(X,P_\omega,\mu).$ We denote by $\P_{\omega,\mu}$ the
corresponding probability measure.

We assume throughout this note that $m^*>1,$ excluding the case
where there is only one particle.

The first result is the following $0$-$1$ law.

\begin{lem}\label{lem:5}
We have either
\begin{itemize}
    \item for $\Theta$-a.a. realizations $(\omega,\mu)$ it holds that
     $\E_{\omega,\mu,x} \nu(x_0)<\infty$ for all $x,x_0\in X,$ or
    \item for $\Theta$-a.a. realizations  $(\omega,\mu)$ it holds that
    $\E_{\omega,\mu,x_0} \nu(x_0)=\infty$ for all $x_0\in X.$
\end{itemize}
\end{lem}
\begin{proof}
Let $x$ and $x_0$ be such that $\E_{\omega,\mu,x}\nu(x_0)=\infty.$
Hence, $\E_{\omega,\mu,x_0}\nu(x_0)=\infty.$ Let $k$ be such that
$\E_{\omega,\mu,x_0} \eta(k,x_0)>1.$ Following the proof of (ii)
in Theorem \ref{thm:1} we obtain that
$\E_{\omega,\mu,y_0}\nu(y_0)=\infty$ for all $y_0$ s.t.
$d(x_0,y_0)>k.$  Assume the distributions $Q$ and $\widetilde{Q}$
to be discrete. Then for $\Theta$-a.a. realizations $(\omega,\mu)$
there exists $z_0,$ $d(z_0,x_0)>2k,$ such that
$p_\omega(z_0x_i,z_0 x_j)=p_\omega(x_i, x_j)$ and
$\mu(z_0x_k)=\mu(x_k)$ for all $x_i,x_j,x_k\in\{x:~d(x_0,x)\leq
k\}.$ Therefore,
$$\E_{\omega,\mu,z_0}\eta(k,z_0)=\E_{\omega,\mu,x_0}
\eta(k,x_0)>1$$ and we conclude  that
$\E_{\omega,\mu,y_0}\nu(y_0)=\infty$ for all $y_0$ s.t.
$d(y_0,x_0)\leq k.$ In order to show the general case one combines
continuity arguments similar to those in the proof of Lemma
\ref{lem:2} with approximations arguments of the type used in the
proof of \ref{thm:1}, part(ii).
\end{proof}
In particular, the BRWRE is either transient for $\Theta$-a.a.
environments or recurrent for $\Theta$-a.a. environments. We have
even the stronger result: \pagebreak
\begin{thm}\label{thm:0-1}
We have either
\begin{itemize}
    \item for $\Theta$-a.a. realizations $(\omega,\mu)$ the BRWRE is strongly
    recurrent:
    $$ \P_{\omega,\mu,x}\left( \sum_{n=1}^\infty
\sum_{i=1}^{\eta(n)}\mathbf{1}\{x_i(n)=x\}=\infty\right)=1\quad
\forall x\in X,\mbox{ or}$$
    \item for $\Theta$-a.a. realizations $(\omega,\mu)$ the BRWRE is transient:
    $$ \P_{\omega,\mu,x}\left( \sum_{n=1}^\infty
\sum_{i=1}^{\eta(n)}\mathbf{1}\{x_i(n)=x\}=\infty\right)=0\quad
\forall x\in X.$$
\end{itemize}
\end{thm}
\begin{proof}
Since the proofs of   Propositions  (1.1) and  (1.2)  in
\cite{comets05} carry over to BRWRE on Cayley Graphs  we just give
a brief sketch how the claim follows from Lemma \ref{lem:5} for
$Q$ and $\widetilde{Q}$ discrete. It suffices to show that
$\E_{\omega,\mu,x_0} \nu(x_0)=\infty$ for all $x_0\in X$ implies
$\alpha\equiv 1.$ Let $k$ be as in the proof of Lemma \ref{lem:5}
and consider a distinguished ancestry line.  At any time a
particle in this ancestry line splits up in at least two particles
we start a Galton-Watson process $(\zeta_i)$ defined as in the
proof of Theorem \ref{thm:1}, part (ii). We obtain a sequence of
Galton-Watson processes $(\zeta_i^n),~n\geq 1.$ Due to the choice
of $k,$ the discreteness of $Q\times\widetilde Q$  and the
infiniteness of $X$ we can extract a subsequence $n_l$ such that
$E[\zeta_1^{n_l}]=\eta(k,x_0)>1$ and such that the processes
$(\zeta_i^{n_l})$ are independent. Hence, at least one of these
processes will survive and $\alpha\equiv 1.$
\end{proof}
We give the classification for BRWRE in transience and strong
recurrence.
\begin{thm}\label{thm:5}
 If $m^*\leq  \frac1{\rho} $ then the BRWRE is transient for
$\Theta$-a.a. realizations $(\omega,\mu)$, otherwise it is
strongly recurrent for $\Theta$-a.a. realizations $(\omega,\mu)$.
\end{thm}
\begin{proof}
Let $m^*\leq m_c:=\frac{1}{\rho}.$ For $\Theta$-a.a. realizations
$(\omega,\mu)$ we couple $(X,P_\omega,\mu)$ with
$(X,P_\omega,\tilde\mu)$ with $\widetilde m=m_c.$ Theorem
\ref{thm:1} implies that $(X,P_\omega,\tilde \mu)$ is transient.
Hence $(X,P_\omega,\mu)$ is transient due to Remark
\ref{rem:couple}.

We shall prove the converse by showing that for $\Theta$-a.a.
realizations $(\omega,\mu)$ there exists no $f>0$ satisfying
$P_\omega f(x) \leq \frac {f(x)}{m(x)}\quad \forall x\in X$ and
$f(x_0)=1.$ We conclude with Theorems \ref{thm:1} and
\ref{thm:0-1}. Note that if $f$ is a solution of $P_\omega f(x)
\leq \frac {f(x)}{m(x)}\quad \forall x\in X$ then $c\cdot f$ is a
solution for all $c\in\R^+$ and we can assume $f(x_0)=1.$

Let $\widetilde{m}:=m^*-\eps,$ with $\eps>0$ such that
$\widetilde{m}>m_c.$
\newline {\bf Claim:} For $\Theta$-a.a. $\omega$ there exists
$K=K(x_0,\omega)>0$ such that for all functions $h$ with $P_\omega
h(x)\leq  \frac{h(x)}{\widetilde{m}}~\forall x$ and $h(x_0)=1$,
$h(x)\leq 0$ for some $x$ with $d(x_0,x)<K,$ where
$d(\cdot,\cdot)$ is the usual graph distance.
\newline {\bf Proof:} Assuming the opposite, we have a
sequence of functions $h_n$ with $h_n(x_0)=1,$

$$P_\omega h_n(x)\leq \frac{h_n(x)}{\widetilde{m}}~\forall x\in X $$
and $h_n(x)>0$ for all $x$ with $d(x_0,x)\leq n.$ Let
$g(x):=\liminf_{n\rightarrow\infty} h_n(x),~ \forall x\in X .$ We
have $g(x_0)=1,$ $g(x)\geq 0$ and with Fatou's Lemma:
\begin{equation}\label{eq:thm:5.1}
P_\omega g(x)\leq \frac{g(x)}{\widetilde{m}},\quad\forall x\in X .
\end{equation}

Since $g(x) \geq {\widetilde{m}}\cdot P_\omega g(x)=
{\widetilde{m}}\cdot \sum_y p_\omega(x,y) g(y)~\forall x,~
g(x_0)=1$ and $P_\omega$ is irreducible we have that $g(x)>0$ for
all $x.$ Equation (\ref{eq:thm:5.1}) together with Theorem
\ref{thm:1} implies that the BMC with constant mean offspring
$\widetilde{m}$ and underlying Markov Chain $(X,P_\omega)$ is
transient for $\Theta$-a.a. $\omega.$ Since
$\widetilde{m}>m_c=1/\rho=1/{\rho(P_\omega)}$ this contradicts
Theorem \ref{thm:1a}. This proves the claim.

We use the independence of $Q$ and $\widetilde Q.$ Let $\omega$ be
a typical realization of the environment. With positive
$\Theta$-probability  the branching rates in $\{y:~ d(x_0,y)\leq
K(x_0,\omega)\}$ are higher than $\widetilde m.$  In this case, we
couple the BRWRE $(X,P_\omega,\mu)$  with a process
$(X,P_\omega,\tilde\mu)$  with fewer particles and mean offspring
$\widetilde m=\widetilde m(x)$ for all $x\in \{y:~ d(x_0,y)\leq
K\}$ and $\widetilde\mu(x)=\mu(x)$ for all $x\in \{y:~ d(x_0,y)>
K\}.$ Due to the Claim there exists no positive function $f$ such
that $P_\omega f(x)\leq \frac{f(x)}{\widetilde m}$ for all $x\in
\{y:~ d(x_0,y)\leq K\}.$ Therefore there exists no positive $f$
such that $P_\omega f(x)\leq \frac{f(x)}{\widetilde m(x)}$ for all
$x\in X.$ Due to Theorem \ref{thm:1} we have recurrence of the
coupled process $(X,P_\omega,\tilde\mu).$ The recurrence of
$(X,P_\omega,\mu)$ follows with Remark \ref{rem:couple}.
Eventually due to Theorem \ref{thm:0-1} we have that BRWRE is
strongly recurrent for $\Theta$-a.a. environments.
\end{proof}

 The transience resp. recurrence does only depend on the support of the
environment. Thus, suppose that a BRWRE is recurrent for almost
all realizations for a marginal distribution $Q_1\times \widetilde
Q_1.$ Then every BRWRE with distribution $Q_2\times \widetilde
Q_2$ such that $supp(Q_1)= supp(Q_2)$ and $supp(\widetilde Q_1)=
supp(\widetilde Q_2)$ is recurrent for a.a. realizations. For
BRWRE on $\Z^d,$ this was already shown in \cite{comets05}.

Furthermore, Theorem \ref{thm:5} states  that the condition
$supp(\widetilde Q_1)= supp(\widetilde Q_2)$ can be replaced by
$m^*_1=m^*_2.$ In the following section we show that we can
replace the condition $supp(Q_1)= supp(Q_2)$ by $conv(supp(Q_1))=
conv(supp(Q_2)).$ Thus,  recurrence and transience only depends on
some extremal points of the support of $Q\times \widetilde Q.$
Varadhan showed in \cite{varadhan2003} that the spectral radius of
a RWRE on $\Z^d$ only depends on the convex hull of the support.
His arguments
 immediately apply to RWRE on Cayley Graphs. We give a modified proof,
 which uses properties of the BRWRE instead of approximations of
 the spectral radius.

\begin{thm}\label{thm:6}
We have for $ \eta$-a.a. $\omega,$
\begin{eqnarray*}
  \rho=\rho(P_{\omega}) &=& \sup_{\sigma} \rho(P_{\sigma}) \\
   &=& \sup_{\hat{\sigma}} \rho(P_{\hat{\sigma}}) ,
\end{eqnarray*}
where the latter $\sup$ is over all possible realizations
${\hat{\sigma}}=({\hat{\sigma}}_x)_{x\in X}$ with
${\hat{\sigma}}_x\in\hat{\mathcal{K}}.$
\end{thm}

\begin{proof}
In order to prove the claim, that is stated for RWRE, we consider
the corresponding BRWRE with a given offspring distribution $\mu$,
$\mu(x)=\mu$ for all $x,$ and mean offspring $m.$  The fact that
the behavior of the BRWRE depends on the mean offspring $m$ is
used frequently. Let $x_0$ be the origin of the corresponding
BMC*.

 For $\eta$-a.a. $\omega$ we have that $m\leq 1/\rho(P_\omega)$ implies
the transience of $(X,P_\omega,\mu)$ and therefore
$\E_{\omega,\mu,x} \nu(x_0)<\infty$ for all $x\in X,$ see Theorem
\ref{thm:1}. In this case the function
$f_{\omega,m}(x):=\E_{\omega,\mu,x} \nu(x_0)$ for $x\neq x_0$ and
$f_{\omega,m}(x_0):=1$ is a solution of
\begin{equation}\label{eq}
 P_\omega f(x)=\frac{f(x)}{m(x)}~\forall x\neq
x_0,
\end{equation}
see Theorem \ref{thm:4}. On the other hand we have for $\eta$-a.a.
$\omega,$ that if there exists a function satisfying equation
(\ref{eq}) we have that $\E_{\omega,\mu,x}\nu(x_0)<\infty,$ for
all $x\neq x_0,$ due to  Theorem \ref{thm:4}. Since $S$ is finite
we obtain that $$ \E_{\omega,\mu,x_0}\nu(x_0)=m\cdot
\sum_{y:~x^{-1}y\in S} p(x,y) \E_{\omega,\mu,y}\nu(x_0)<\infty,$$
too. Due to Lemma \ref{lem:5} and Theorem \ref{thm:1} this implies
the transience of the BMC $(X,P_\omega,\mu)$ and eventually that
$m\leq 1/\rho.$ Therefore, for $\eta$-a.a. $\omega$ the existence
of a function $f$ satisfying equation (\ref{eq}) is equivalent to
the transience of $(X,P_\omega,\mu).$ Note that this equivalent to
$m\leq 1/\rho.$ To make use of this fact we investigate the values
of $m$ such that $f_{\sigma,m}(x):=\E_{\sigma,\mu,x} \nu(x_0)$  is
finite for all possible realizations $\sigma.$ We show that there
exists a critical $\widetilde m$ such that $f_{\sigma,m}(x)$ is
finite for all $\sigma$ if $m<\widetilde m$ and infinite for some
$\sigma$ if $m>\widetilde m.$ Using the properties of BRWRE we
show  that $\widetilde m= 1/\rho.$ The claim will then follow by
considering the BRWRE with support $\hat{\mathcal{K}}$ and
observing that $f(\hat\sigma,x)$ is finite for all
$\hat\sigma=(\hat\sigma_x)_{x\in X}$ if $m<\widetilde m$ and
infinite for some $\hat\sigma$ if $m>\widetilde m.$ Let us first
show that there exists a critical $\widetilde m.$

Let $\sigma$ be a possible realization
${\sigma}=({\sigma}_x)_{x\in X}$ with ${\sigma}_x\in\mathcal{K}$
and consider the equation
\begin{eqnarray}\label{eq:star}
  P_\sigma f(x) &=& \sum_y p_\sigma(x,y) f(y) =\frac{f(x)}m,
  \quad \forall x\neq x_0 \mbox{ and } f(x_0)=1.
\end{eqnarray}
If $\E_{\sigma,\mu,x}\nu(x_0)<\infty$ for $x\neq x_0$ we can write
a solution of (\ref{eq:star}) as
$$f_{\sigma,m}( x)=\E_{\sigma,\mu,x}\nu(x_0)\mbox{ for }x\neq x_0
\mbox{ and } f_{\sigma,m}(x_0)=1,$$ due to Theorem \ref{thm:4}.
 In the
following, we consider ${\sigma}=({\sigma}_x)_{x\in X}$ with
${\sigma}_x\in\mathcal{K}$ as a choice, chosen at will, of
transition probability functions and $f_{\sigma,m}(x)$ as a payoff
function to be maximized. To show the finiteness of
$f_{\sigma,m}(x)$ for all $\sigma$ we maximize the function
$f_{\sigma,m}(x)$ in $\sigma$ and determine those values of $m$
such that this maximum is finite. This is a typical problem of
dynamic programming. The corresponding Bellman equation is:
\begin{eqnarray}\label{eq:2star}
  m\cdot \sup_{p(x,\cdot)\in \mathcal{K}} \sum_y p(x,y) f(y) &=& f(x)
  \quad \forall x\ne x_0 \\
 \nonumber f(x_0)&=&1
\end{eqnarray}
This problem can be understood as an infinite stage allocation
process. Since the existence of a function satisfying
(\ref{eq:2star}) is  not  guaranteed, we  first consider the
$N$-stage allocation process:
\begin{eqnarray}
\nonumber
  f_n(x) &=& m\cdot \sup_{p(x,\cdot)\in \mathcal{K}}
  \sum_y p(x,y) f_{n-1}(y) \quad \forall x\ne x_0\quad 1<n\leq N \\
 \nonumber  f_n(x_0)&=&f_{n-1}(x_0)\quad 1<n\leq N\\
 \nonumber  f_1(x) &=& \delta_{x_0}
\end{eqnarray}
Observe that the sequence $\{f_N(x)\}_{N\in\N}$ is  increasing for
all $x\in X$. Hence, there exists a largest
$\widetilde{m}\in\R\cup\{\infty\}$ such that $f_N(x)$ is bounded
for all $m<\widetilde{m}$ and hence due to the monotone
convergence theorem we have that $f(x):=\lim_{N\rightarrow \infty}
f_N(x)$ exists and verifies (\ref{eq:2star}) for all
$m<\widetilde{m}.$ The latter can be shown via a standard
argument: We have by monotonicity in $N$ and
$$f_{N+1}(x)=m\cdot \sup_{p(x,\cdot)\in \mathcal{K}} \sum_y p(x,y)
f_N(y)$$that for all $N\in\N:$
\begin{eqnarray*}
  f(x) &\geq & m \cdot \sup_{p(x,\cdot)\in \mathcal{K}} \sum_y p(x,y) f_N(y)  \\
   &\geq& m \cdot \sum_y p(x,y) f_N(y)\quad \forall p(x,\cdot)\in\mathcal{K}
\end{eqnarray*}
Letting $N\rightarrow\infty$, this yields
$$f(x) \geq  m \cdot \sum_y p(x,y)f(y)\quad \forall
p(x,\cdot)\in\mathcal{K},$$

and hence $$f(x)\geq m\cdot\sup_{p(x,\cdot)\in \mathcal{K}} \sum_y
p(x,y)f(y).$$ On the other hand we have for all $N:$
\begin{eqnarray}
  f_{N+1}(x) &=& m\cdot \sup_{p(x,\cdot)\in \mathcal{K}} \sum_y p(x,y) f_N(x) \\
   &\leq& m\cdot \sup_{p(x,\cdot)\in \mathcal{K}}
    \sum_y p(x,y)f(y)
\end{eqnarray}
and therefore $f(x)\leq m\cdot \sup \sum_y p(x,y)f(y).$ It remains
to show that $\widetilde{m}=1/\rho.$ Let $\widetilde{f}_m$ be the
solution of the Bellman equation (\ref{eq:2star}) dependent on the
parameter $m$. For each $m$ and $\sigma$, we consider the
corresponding BMC $(X,P_\sigma,\mu)$ with constant mean offspring
$m$ and transition probabilities $P_\sigma.$ Let
$m<\widetilde{m},$ hence $\widetilde{f}_m$ exists and is finite
and so does $f_{\omega,m}(x)=\E_{\omega,\mu,x}\nu(x_0)$ for
$\eta$-a.a. realizations $\omega,$ since $\widetilde{f}_m$ is
maximal. Due to Lemma \ref{lem:5} and Theorem \ref{thm:1} we have
that the corresponding BMC $(X,P_\omega, \mu)$ is transient for
$\eta$-a.a. $\omega.$ Thus, with Theorem \ref{thm:5},
$\widetilde{m}\leq 1/\rho.$ In order to show the converse observe
that $m\leq 1/\rho$ implies that $m\leq 1/\rho(P_\sigma)$ for all
$\sigma,$ see Lemma \ref{lem:2}. Hence, for $m\leq 1/\rho$  the
BMC $(X,P_\sigma,\mu)$ is transient for all possible realizations
$\sigma,$ see Theorem \ref{thm:5}. Due to Theorem \ref{thm:1} we
have $f_{\sigma,m}(x)=\E_{\sigma,\mu,x}\nu(x_0)<\infty$ for all
$x\in X$ and all $\sigma=(\sigma_x)_{x\in X}$ and hence
$\widetilde m\geq 1/\rho.$

It is obvious that the value $\widetilde{m}$ and hence the value
$\rho$ does not change if we consider the BRWRE with support
$\hat{\mathcal{K}}$ and, instead of equation (\ref{eq:2star}), the
dynamic programming problem:
\begin{eqnarray}\nonumber
  m\cdot \sup_{p(x,\cdot)\in \hat{\mathcal{K}}} \sum_y p(x,y) f(y) &=& f(x)
  \quad \forall x\ne x_0 \\
 \nonumber f(x_0)&=&1,
\end{eqnarray} where the $\sup$ is over the convex hull $\hat{\mathcal{K}}$  of
$\mathcal{K}.$
\end{proof}

\subsection{BRWRE on $\Z^d$}
We consider the case where $G=\Z^d$ and $S$ is some finite
generator of $\Z^d.$ Thus  the process becomes the BRWRE on $\Z^d$
with bounded jumps. In this case one can explicitly calculate the
spectral radius $\rho.$ We follow the argumentation of
\cite{varadhan2003} to show the following Lemma:

\begin{lem}\label{lem:6}
For a RWRE on $\Z^d$ we have for $\eta$-a.a. realizations $\omega$
\begin{eqnarray}
   \rho(P_\omega) &=& \sup_{p\in\hat{\mathcal{K}}} \rho(P^h_p) \\
   &=& \sup_{p\in\hat{\mathcal{K}}} \inf_{\theta\in\R^d}
   \left( \sum_{s\in S} e^{\langle
   \theta, s\rangle} p(s)\right),
\end{eqnarray}
where $P^h_p$ is the transition matrix of the random walk with
$p(x,x+s)=p(0,s)=:p(s)$ for all $x\in\Z^d,~s\in S.$ In particular,
we have that $\rho=1$ if and only if there is a
$p\in\hat{\mathcal{K}}$ with $\sum_s s p(s)=0.$
\end{lem}
\begin{proof}
The second equality is more or less standard. It  follows for
example  from the fact that $\rho(P)=\exp(-I(0)),$ where
$I(\cdot)$ is the rate function of the large deviations of the
random walk determined by $P.$ One direction of the first equality
follows directly from the Lemma \ref{lem:1} and Theorem
\ref{thm:6}. Thus it remains to show that
$$\rho(P_\omega)\leq  \sup_{p\in\hat{\mathcal{K}}}
\inf_{\theta\in\R^d} \left(\sum_s e^{\langle
   \theta, s\rangle} p(s)\right)$$ for $\eta$-a.a. realizations $\omega.$
Observing that the function $\phi(p(\cdot),\theta):=\left(\sum_s
e^{\langle \theta, s\rangle} p(s)\right)$ is convex in $\theta$
and linear in $p(\cdot)$, we get by a standard minimax argument
that
$$\sup_{p\in\hat{\mathcal{K}}} \inf_{\theta\in\R^d} \sum_s e^{\langle
   \theta, s\rangle} p(s)= \inf_{\theta\in\R^d} \sup_{p\in\hat{\mathcal{K}}} \sum_s e^{\langle
   \theta, s\rangle} p(s)=:c.$$ Let $\eps>0$ and $\theta\in \R^d$
   such that
   $$ \sup_{p\in\hat{\mathcal{K}}} \sum_s e^{\langle
   \theta, s\rangle} p(s) \leq c (1+\eps). $$ By induction we have
   for any realization $\omega:$
   $$ E_\omega\left[\exp(\langle \theta, X_n\rangle)  \right]\leq
   (c (1+\eps))^n,$$ where $X_n$ is the position of the RW at time $n.$ Therefore by observing the event  $\{X_n=0\}:$
    $$ p_\omega^{(n)}(0,0)\leq (c (1+\eps)) ^n,$$ and hence
    $\rho(P_\omega)\leq c (1+\eps)$ for all $\eps>0.$

The last part of the Lemma follows now from the observation that
$\inf_{\theta\in\R^d} \sum_s e^{\langle
   \theta, s\rangle} p(s)=1$ if and only if $\sum_s s p(s)=0.$
\end{proof}

We immediately obtain the following criteria.
\begin{cor}\label{cor:6}
The BRWRE is strongly recurrent for $\Theta$-a.a. realizations if
$$ (m^*)^{-1} <\sup_{p\in\hat{\mathcal{K}}} \inf_{\theta\in\R^d} \left( \sum_s e^{\langle
   \theta, s\rangle} p(s)\right).$$ Otherwise it is transient for
   $\Theta$-a.a. realizations.
\end{cor}

\begin{ex}
We consider the nearest-neighbor BRWRE on $\Z^d.$ Let $S=\{\pm
e_i,~1\leq i\leq d\}.$ The BRWRE is strongly recurrent for
$\Theta$-a.a. realizations $(\omega,\mu)$ if
$$ (m^*)^{-1} < \sup_{p\in\hat{\mathcal{K}}}
 \left(2 \sum_{i=1}^d \sqrt{p(e_i) p(-e_i)}  \right).$$
Otherwise it is transient for
   $\Theta$-a.a. realizations.
\end{ex}

\subsubsection*{Acknowledgment}
The author is grateful to  Nina Gantert for valuable discussions
and her helpful comments on a previous version of this paper.

\bigskip
 \noindent
\begin{tabular}{l}
Sebastian M\"uller \\
 Institut f\"ur Mathematische Statistik\\
 Universit\"at Münster\\
 Einsteinstr. 62  \\
 D-48149 Münster\\
 Germany\\
{\tt Sebastian.Mueller@math.uni-muenster.de}\\
\end{tabular}

\end{document}